\documentclass[10pt]{article}
\usepackage{amsmath,amssymb,graphicx, epsfig, psfrag,color}
\usepackage[boxruled]{algorithm2e}

\newtheorem{thm}{Theorem}[section]
\newtheorem{lem}[thm]{Lemma}
\newtheorem{cor}[thm]{Corollary}
\newtheorem{defn}[thm]{Definition}

\newcommand{\be}{\begin{enumerate}}
\newcommand{\ee}{\end{enumerate}}
\newcommand{\bc}{\begin{center}}
\newcommand{\ec}{\end{center}}
\newcommand{\bt}{\begin{tabular}}
\newcommand{\et}{\end{tabular}}
\newcommand{\ba}{\begin{array}}
\newcommand{\ea}{\end{array}}

\newcommand{\noi}{\noindent}
\newcommand{\Z}{\mathbb Z}

\newcommand{\bp}{\noindent\textit{Proof:}  }
\newcommand{\ep}{\hfill$\Box$}

\newcommand{\MDC}{\textrm{\textsc{MDC}}}
\def\Area{\hbox{\rm Area}}
\def\onto{{\kern3pt\to\kern-8pt\to\kern3pt}}
\newcommand{\e}{\varepsilon}
\newcommand{\hnn}{HNN-extension}
\def\ss{\smallskip}
\def\ms{\medskip}
\def\bs{\bigskip}

\DeclareMathOperator{\hght}{height}

\definecolor{dmagenta}{rgb}{.5,0,.5} 
\definecolor{dred}{rgb}{.5,0,0} 
\definecolor{dgreen}{rgb}{0,.5,0} 
\definecolor{dblue}{rgb}{0,0,0.5} 
\definecolor{black}{rgb}{0,0,0} 
\definecolor{vdgreen}{rgb}{0,.3,0} 
\definecolor{vdred}{rgb}{.3,0,0} 
\definecolor{red}{rgb}{1,0,0} 

\makeatletter
\newcounter{algorithmbis}
\renewcommand{\thealgorithmbis}{\arabic{algorithmbis}}
\def\algorithmbis{\@ifnextchar[{\@algorithmbisa}{\@algorithmbisb}}
\def\@algorithmbisa[#1]{%
 \refstepcounter{algorithmbis}
 \trivlist
 \leftmargin\z@
 \itemindent\z@
 \labelsep\z@
 \item[\parbox{\textwidth}{%
   \hrule
   \hrule
   \noindent\strut\textbf{Algorithm \thealgorithmbis \ } #1      
   \hrule
 }]\hfil\vskip0em%
}
\def\@algorithmbisb{\@algorithmbisa[]}

\makeatother

\SetKwInput{KwIn}{Input}
\SetKwInput{KwConv}{Method\rule{0mm}{7mm}}  
\SetKwInput{KwOut}{Goal\rule{0mm}{7mm}}  

\begin{document}
\title{The Dehn function of Stallings' group \\ 
\rule{0mm}{8mm}\small{To appear in \emph{Geometric and Functional Analysis}}}

\author
{Will Dison\\
\small Department of Mathematics, \\[-0.8ex]
\small The University of Bristol, University Walk, Bristol, BS8 1TW, United Kingdom  \\[-0.8ex]
\small\texttt{w.dison@bristol.ac.uk}\\[1.2ex]
Murray Elder\\
\small Mathematics,\\[-0.8ex]
\small The University of Queensland, Brisbane, Queensland 4072, Australia \\[-0.8ex]
\small\texttt{m.elder@uq.edu.au}\\[1.2ex]
Timothy R. Riley\footnote{The third author is grateful for support
from NSF grant DMS--0540830 and for the hospitality of the Institut des Hautes \'Etudes
Scientifiques in Paris during the writing of this article.}
\footnote{Corresponding author}\\
\small Department of Mathematics, \\[-0.8ex]
\small  The University of Bristol, University Walk, Bristol, BS8 1TW, United Kingdom \\[-0.8ex]
\small\texttt{tim.riley@bris.ac.uk} \\[1.2ex]
Robert Young\\
\small Institut des Hautes \'Etudes Scientifiques,\\[-0.8ex]
\small Le Bois Marie, 35 route de Chartres, F-91440 Bures-sur-Yvette, France \\[-0.8ex]
\small\texttt{rjyoung@ihes.fr}
}

\date{\small
}

\maketitle
\newlength{\miniwidth} \addtolength{\miniwidth}{\linewidth} \addtolength{\miniwidth}{-\parindent}

\vspace*{-6mm}

\begin{center}
\emph{Dedicated to John Stallings.} \\ 
\end{center}

\begin{abstract}
\noi We prove that the Dehn function of a group of Stallings that is finitely presented but not of type $\mathcal{F}_3$ is quadratic.

\end{abstract}

\section{Introduction}\label{sec:intro}

\noi
A group is of type $\mathcal F_1$ when it can be finitely generated, $\mathcal F_2$ when it can be finitely presented, and more generally $\mathcal F_n$ when it admits an Eilenberg--Maclane space with finite $n$-skeleton.  In the early 1960s Stallings \cite{Stallings} constructed a group $S$ that is $\mathcal{F}_2$ but not $\mathcal{F}_3$.
  Bieri~\cite{Bieri} recognised $S$ to be
\begin{equation} \label{kernel}
\textup{Ker}( \, F(\alpha, \beta) \times F(\gamma, \delta) \times F(\epsilon, \zeta) \ \onto \ \Z \, )
\end{equation}
where the map is that from the product of three
rank-$2$ free groups to $\Z = \langle t \rangle$
which sends all six generators to $t$, and he showed that replacing
$(F_2)^3$ by $(F_2)^n$ gives a family of groups (the
\emph{Bieri--Stallings} groups) of type $\mathcal{F}_{n-1}$ but not
$\mathcal{F}_n$ \cite{Bieri}.

Isoperimetric functions (defined below) for $S$ have been established by a number of authors.
Gersten proved that for $n \geq 3$, the Bieri--Stallings groups admit quintic
isoperimetric functions \cite{Gersten5}; this was sharpened to cubic by Baumslag, Bridson, Miller \& Short in the case of $S$ \cite[\S6]{BBMS}.
Bridson \cite{Bridson4} showed that the
Bieri-Stallings groups were examples of a construction called
\emph{doubling} and argued that a class of doubles including these
groups should also have quadratic isoperimetric functions.
But Groves found an error in his proof \cite{BridsonPersonal, GrovesPersonal}, and it seems that  Bridson's approach, in fact, gives a cubic isoperimetric function, generalising the result in \cite{BBMS}. In this article we establish a quadratic isoperimetric function for $S$, and as $S$ is not hyperbolic (as it is not of type $\mathcal F_3$, for example) this is best possible.  And so we prove:

\begin{thm} \label{Main Theorem}
The Dehn function of Stallings' group is quadratic.
\end{thm}

More precisely, this theorem says that the Dehn function (defined below) of any finite presentation of Stallings' group is equivalent to $n \mapsto n^2$ in the following sense.  For $f,g : \mathbb{N} \to  \mathbb{N}$, we write $f \preceq g$ when $\exists C>0,  \forall n \in \mathbb{N}, f(n) \leq C g(Cn+C) + C n +C$, and we write $f \simeq g$ when $f \preceq g$ and $g \preceq f$.  As is well-known, any two finite presentations of the same group have equivalent Dehn functions.

Our theorem  fulfils Bridson's aim in \cite{Bridson4} of exhibiting wild behaviour within the class of groups with quadratic Dehn functions ---

\begin{cor}
There exists a group with quadratic Dehn functions that is not of type $\mathcal{F}_3$.
\end{cor}

Combined with results in \cite{Papasoglu, Riley3}, the theorem also has the following corollaries.

\begin{cor}
The asymptotic cones of Stallings' group are all simply connected, but not all are 2-connected.
\end{cor}

\begin{cor}  Stallings' group admits a linear isodiametric function.    Indeed,
its filling length function is linear \textup{(}that is, equivalent to $n \mapsto n$\textup{)}.
\end{cor}

 We will work with the presentation
\begin{equation} \label{pres}
 \langle \; a,b,c,d,s \; \mid \; [a,c], \; [a,d], \; [b,c], \; [b,d], \; s^a=s^b=s^c=s^d  \; \rangle
\end{equation}
for $S$ of \cite{BBMS, Gersten5}.  Our notation is $[x,y]:=x^{-1}y^{-1}xy$,  $x^y:= y^{-1}xy$, $x^{-y}:= y^{-1}x^{-1}y$, and  $s^a=s^b=s^c=s^d$  is shorthand for the six defining relations $s^as^{-b}$,  $s^as^{-c}$, $s^as^{-d}$, $s^bs^{-c}$,  $s^bs^{-d}$,  $s^cs^{-d}$. Note that these six relations can be rewritten as $[s,ab^{-1}], [s,ac^{-1}]$, and so on.
  One can view $S$ as an \hnn\ of the product of free groups $F(a,b) \times F(c,d)$ with stable letter $s$ commuting with all elements represented by words on $a^{\pm 1}, b^{\pm 1}, c^{\pm 1}, d^{\pm 1}$ of zero exponent-sum.  The first four relations in the presentation are then the relations coming from $F(a,b)\times F(c,d)$, which we call \emph{commutator relations}.
[Gersten~\cite{Gersten5} showed that this is related to the description of $S$ as a kernel (\ref{kernel}) via $a= \epsilon \alpha ^{-1},  b= \epsilon \beta ^{-1},  c =  \epsilon \gamma^{-1},  d= \epsilon \delta^{-1}, s= \zeta \epsilon^{-1}$.]

We prove Theorem~\ref{Main Theorem} by presenting  an algorithm
(Algorithm~\ref{alg:main}) which takes as input a null-homotopic
word of length $n$ and transforms it to the empty word $\e$ by applying relations from the
presentation~(\ref{pres}). We
call the number of relations applied the \emph{cost} and we wish to
design the algorithm so that this is bounded by a constant multiple of
$n^2$.

To understand the structure of the algorithm, it helps to understand
the structure of a word which represents the identity. Since $S$ is
an HNN-extension of $F(a,b)\times F(c,d)$ by a generator $s$, by
Britton's Lemma, a word $w$ representing the identity contains
``pinches'', or pairs of letters $s$ and $s^{-1}$ separated by a
word which commutes with $s$. (We will later call such words
\emph{balanced}.) Reducing $w$ to the identity involves removing
these pinches by bringing $s$'s and $s^{-1}$'s together.

Using the presentation, one can show that $s$ commutes with words of
the form $xy^{-1}$, where $x,y\in \{a,b,c,d\}$; we will call a
product of such words an \emph{alternating word}.  Then $s$ can
easily be commuted past a product of such words.  The basic strategy
of the algorithm is to identify a pinch, convert the balanced word
inside to an alternating word, and cancel an $s$ and $s^{-1}$.  If
there is a larger pinch containing the alternating word, we can
repeat the process.  Once we have removed all occurrences of the
letters $s$ and $s^{-1}$ from $w$ the resulting word will represent
the identity in $F(a,b) \times F(c, d)$ and we will be able to apply
commutator relations to convert this to the empty word.  Provided
that the process up to this point has not increased the length of
the word significantly, the cost of this final step will be
proportional to $n^2$.

The step which has the highest cost is converting a balanced word to
an alternating word.  One way to do this involves first separating
$a$'s and $b$'s from $c$'s and $d$'s, then inserting $a$'s and $c$'s
to produce an alternating word.  For a word in $a,b,c,d$ of length
$l$, this has cost approximately $l^2$. If the pinches are deeply
nested, we will need to repeat the process up to $n/2$ times, and as
the second step largely undoes the first the total cost could be up
to $n^3$.

We improve this by employing two key techniques.  First, we utilise
a divide-and-conquer strategy to convert balanced words to
alternating form.  We partition a balanced word into subwords and
separate $a$'s and $b$'s from $c$'s and $d$'s in each subword,
rather than in the whole word, before inserting $a$'s and $c$'s to make it alternating.  We will say that the resulting word
is in \emph{partitioned alternating form}.


A typical intermediate stage in our process is a word with several
subwords in partitioned alternating form.  Indeed, we will specify
intermediate stages by a list of subwords of the original word and
partitions of these subwords; the intermediate stage will then be the
original word with the specified subwords replaced by their
partitioned alternating forms.  As the algorithm progresses, these
subwords grow and we merge adjacent subwords and adjacent pieces of
the partitions.  When two pieces in a partition of a balanced subword
are merged, the cost is proportional to the square of the length of
the words.

\textit{A priori}, these merges could have a heavy total cost.  To overcome
this problem we employ a second key technique: we only use a
particular type of partitioned alternating form, which we call
\emph{dyadic alternating form}.  In this form, the partition
of a balanced word only involves subwords of length $2^k$, where $k
\in \mathbb{N}$.  Since once two pieces are merged together, they are never separated, there can be at most $n/2^k$ merges of pieces of length $2^{k-1}$, and thus the total cost of all mergings will be
proportional to
$$
\sum_{k=1}^{\lceil\log n\rceil} \frac{n}{2^k} (2^{k-1})^2 \ \simeq \ n^2.$$

This article is organised as follows. In Section~\ref{sec:prelim} we define alternating and balanced words, and we establish some basic facts about them. In Section~\ref{sec:dyadic} we define \emph{dyadic alternating form} --- this involves breaking up balanced words using a \emph{dyadic partition}.   
Our main algorithm is in Section~\ref{sec:main} and is analysed in  Section~\ref{sec:cost}.  It proceeds by converting more and more of the input word $w$ into dyadic alternating form by calling a number of subroutines (given in Section~\ref{sec:subroutines}) to combine smaller subwords in dyadic alternating form into larger ones.

\bs \noi \emph{Article history}. A number of prior versions of this article were made public.  In the first, Elder and Riley established that  $n^{5/2}$ is an isoperimetric function for $S$.  Dison realised the result could be improved to $n^{7/3}$ and produced a new version of the paper in collaboration with Elder and Riley.  Later Young contributed further insights that achieve the definitive $n^2$ result, and he, together with the other three authors, produced this version.

\bs \noi \emph{Acknowledgements.}   We are grateful to Noel~Brady, Martin~Bridson,
Daniel~Groves and Steve~Pride for discussions on this problem, and to an anonymous referee for a careful reading.

\section{Preliminaries} \label{sec:prelim}

 \noi   Write $u=u(a_1, \ldots, a_k)$ when $u$ is a word on ${a_1}^{\pm 1}, \ldots, {a_k}^{\pm 1}$.  Write $u[i]$ for the $i$-\text{th} letter of $u$.
The length of $u$ as a word (with no free reductions performed) is $\ell(u)$.  
The sum of the exponents of the letters in $u$ is denoted by  $\xi(u)$, which we call the {\em exponent sum} of $u$.  Unless otherwise indicated, we consider two words to be equal when they are identical letter-by-letter.  A \emph{partition} of $u$ is any way of expressing $u$ as a  concatenation  $u_1 \ldots u_k$ of subwords.   We denote the empty word by $\e$.

\begin{defn}[Cost, Dehn function, isoperimetric function] \label{cost etc}
Given words $w$, $w'$ representing the same element of a group with
finite presentation \mbox{$\langle \mathcal{A} \mid \mathcal{R}
\rangle$}, one can convert $w$ to $w'$ via a sequence
of words $W=(w_i)_{i=0}^m$ in which $w_0=w$, $w_m=w'$ and for each
$i$, $w_{i+1}$ is obtained from $w_i$ by free reduction
\textup{(}$w_i = \alpha a a^{-1} \beta \mapsto \alpha \beta
=w_{i+1}$ where $a \in \mathcal{A}^{\pm 1}$\textup{)}, by free
expansion \textup{(}the inverse of  a free reduction\textup{)}, or
by applying a relator \textup{(}$w_i=\alpha u\beta \mapsto \alpha
v\beta = w_{i+1}$ where a cyclic conjugate of $uv^{-1}$  is in
$\mathcal{R}^{\pm 1}$\textup{)}.   The \emph{cost} of $W$ is the
number of $i$ such that $w_i \mapsto w_{i+1}$ is an
\emph{application-of-a-relator} move.  For words $w$ that represent
the identity \textup{(}i.e.\ \emph{null-homotopic} words\textup{)},
 $\Area(w)$ is  the minimal cost amongst all $W$ converting $w$ to $\e$.  The
\emph{Dehn function} $\delta : \mathbb{N} \to \mathbb{N}$ of $\langle \mathcal{A} \mid \mathcal{R} \rangle$ is
$$\delta(n) \ := \  \max \{ \Area(w) \mid w  \textup{ represents } 1 \textup{ and } \ell(w) \leq n\}.$$
An \emph{isoperimetric function} for $\langle \mathcal{A} \mid \mathcal{R} \rangle$ is any $f:  \mathbb{N} \to \mathbb{N}$  such that  $\delta(n) \leq   f(n$) for all $n$.
\end{defn}

\begin{defn}[Alternating words]\label{alternating}
A word $u = u(a,b,c,d)$ is \emph{alternating} if $u$ has even length and $u[i]$ is in  $\{a,b,c,d\}$ for all odd $i$ and in $\{a^{-1},b^{-1},c^{-1},d^{-1}\}$ for all even $i$.
\end{defn}

\noi [The reader familiar with van~Kampen diagrams and corridors (also known as bands) may find it helpful to note that alternating words are those which, after removing all $aa^{-1}, bb^{-1}, cc^{-1}$ and $dd^{-1}$ subwords, can be read along the sides of  $s$-corridors in van~Kampen diagrams over $S$.]

\begin{defn}[Balanced words]  \label{balanced}
A word $u=u(a,b,c,d,s)$ is \emph{balanced} if there exists an alternating word $v = v(a,b,
c, d)$ with $u = v$ in $S$.
\end{defn}

\begin{lem}
For a word $u=u(a,b,c,d,s)$, the following are equivalent.
\begin{itemize}
\item[\textup{(}\textit{i}\textup{)}]   $u$ is balanced.
\item[\textup{(}\textit{ii}\textup{)}]  $u$ represents an element of $\langle a,b,c,d \rangle$ in $S$ that commutes with $s$.
\item[\textup{(}\textit{iii}\textup{)}]  $u$ represents an element of $\langle a,b,c,d \rangle$ in $S$ and $\xi(u)=0$.
\end{itemize}
\end{lem}

\bp
The equivalence of (\textit{i}) and (\textit{ii}) is straight-forward.

Alternating words have exponent-sum zero so (\textit{i}) implies
(\textit{iii}) by the following observation.  Every relation  of
presentation (\ref{pres}) has exponent-sum zero, so exponent-sum is
preserved whenever a relation is applied to a word and
hence if two words on  $a^{\pm 1},b^{\pm 1},c^{\pm
1},d^{\pm 1}, s^{\pm 1}$ represent the same element in $S$, then
they have the same exponent sum.

To see that  (\textit{iii}) implies (\textit{i}) we can convert a
word $w=w(a, b, c, d)$ with $\xi(w)=0$ to a word in alternating form
as follows. Commute all $a,b$ letters to the front to give a word
$\rho(a,b)\sigma(c,d)$. For each letter of $\rho(a,b)$, replace $a,
b, a^{-1},b^{-1}$ by $ac^{-1}, bc^{-1}, ca^{-1},cb^{-1}$,
respectively. For each letter of $\sigma(c,d)$, replace $c,d,
c^{-1},d^{-1}$ by $ca^{-1}, da^{-1}, ac^{-1},ad^{-1}$, respectively.
In the middle insert $(ca^{-1})^{\xi(\rho(a,b))}$
which cancels out the $c,a$ letters that were added. \ep

\begin{lem} \label{lem:balanced-alternating}
Suppose a word $w= w(a,b,c,d,s)$ is  expressed as $w=\alpha v \beta$
in which $v$ is a balanced subword.  Then $w$  is balanced if and
only if $\alpha\beta$ is balanced.
\end{lem}

\bp Induct on the number of $s$ letters in $w$. The base case where
$w=w(a,b,c,d)$ is immediate
 and the induction step  an application of Britton's Lemma.
Alternatively, this result is  an observation on the layout of
$s$-corridors in a van~Kampen diagram demonstrating that
$w$ equates to some alternating word in $S$. \ep

\begin{lem} \label{lem:xy}
  If $w$ is a balanced word of length at least 2, it contains a subword $xy$ such that either $xy= s^{\pm 1}s^{\mp 1}$ or $x,y \in \{a,b,c,d\}^{\pm 1}$ with $x$ and $y$ having opposite exponents.
\end{lem}
\bp
  By Britton's Lemma, $w$ contains either a subword $s^{\pm 1}s^{\mp
    1}$ or a non-empty balanced subword $u=u(a,b,c,d)$.   In the second case, since $u$ is balanced it has exponent-sum zero and so  must contain a subword $u[i]u[i+1]$ where $u[i]$ and $u[i+1]$ have opposite exponent; take $xy = u[i]u[i+1]$.
    \ep

\section{Dyadic alternating form}\label{sec:dyadic}

The main algorithm will systematically convert subwords of the input word into a special alternating form, which we describe in this section.

\begin{figure}[ht!]
\psfrag{0}{$0$} \psfrag{1}{$1$} \psfrag{2}{$2$} \psfrag{3}{$3$}
\psfrag{4}{$4$} \psfrag{5}{$5$} \psfrag{6}{$6$} \psfrag{7}{$7$}
\psfrag{00}{$D_{0,0}$} \psfrag{01}{$D_{0,1}$}
\psfrag{02}{$D_{0,2}$} \psfrag{03}{$D_{0,3}$}
\psfrag{04}{$D_{0,4}$} \psfrag{05}{$D_{0,5}$}
\psfrag{06}{$D_{0,6}$} \psfrag{07}{$D_{0,7}$}
\psfrag{10}{$D_{1,0}$} \psfrag{11}{$D_{1,1}$}
\psfrag{12}{$D_{1,2}$} \psfrag{13}{$D_{1,3}$}
\psfrag{20}{$D_{2,0}$} \psfrag{21}{$D_{2,1}$}
\psfrag{30}{$D_{3,0}$}
\centerline{\epsfig{file=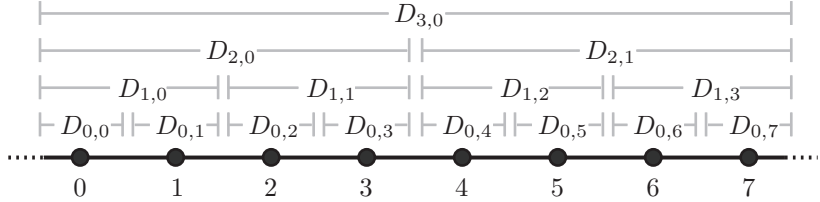}} \caption{Dyadic partitions of the integers.} \label{fig:dyadic}
\end{figure}

By an \emph{interval} in $\Z$, we mean $\Z \cap [\lambda, \mu]$ for some $\lambda, \mu \in \mathbb{R}$.
For $r,j \in \Z$ with $r \geq 0$, define the \emph{dyadic} interval
$D_{r,j}$ to be $\Z \cap [j 2^r ,  (j+1) 2^r)$, as illustrated in Figure~\ref{fig:dyadic}.
Note that $D_{r+1,j}= D_{r,2j}\cup D_{r,2j+1}$ and any two
$D_{r,j}$ and $D_{r',j'}$ are either disjoint or one contains
the other. The {\em height} of $D_{r,j}$ is $r$.

Given an interval $U \subset \Z$, a \emph{cover} of
$U$ is a collection $\{U_i\}_{i=1}^n$ of disjoint intervals $U_i
\subseteq U$ with $\cup_{i=1}^n U_i = U$.  We say that the indexing
on $\{U_i\}_{i=1}^n$ is \emph{ascending} if $\min\left(U_{i+1}\right) = 1 + \max\left(
U_i\right)$ for each $i$.  Given a subword $v$ of a word $w$, let $V
\subset \Z$ be the interval consisting of the positions of the
letters of $v$ in $w$.  Then there is an obvious 1-1 correspondence
between partitions of $v$ and covers of $V$.

\begin{figure}[ht!]
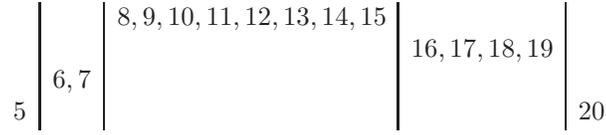

\[\ba{c|c|c|c|c}
& & 8,9,10,11,12,13,14,15 & &  \\
& & & 16,17,18,19  &  \\
& 6,7 & & & \\
5 &    &  &    & 20\\
\ea
\]
\caption{The minimal dyadic cover of $\Z\cap [5,20]$. \label{fig:MDC}}
\end{figure}

A \emph{dyadic cover} of an interval $U \subset \Z$ is a cover
consisting of dyadic intervals.  Define the \emph{minimal dyadic
cover} $\MDC(U)$ of $U$ to be the set of maximal elements (with
respect to containment) of $\{D_{r,j}|D_{r,j}\subseteq U\}$. It is clear that this set is 
 a dyadic cover of $U$. In fact, as a consequence of Lemma~\ref{lem:dyadicToMDC},  it is
 the dyadic cover with the minimal
number of elements.
As an 
example, if $U = \Z\cap [5,20]$, then $\{D_{r,j}|D_{r,j}\subseteq U\}$ is
$$\{D_{0,5},D_{0,6},\ldots, D_{0,20},
D_{1,3},D_{1,4},\ldots, D_{1,9}, D_{2,2},D_{2,3}, D_{2,4}, D_{3,1}
\}$$ and
\begin{eqnarray*}
\MDC(U)  & = &
\{D_{0,5}, D_{1,3},  D_{3,1},  D_{2,4}, D_{0,20}\} \\
& = & \{5\}\cup\{6,7\}\cup\{8,\ldots, 15\}\cup\{16,17,18,19\}\cup\{20\}.
\end{eqnarray*}
 In Figure~\ref{fig:MDC} we display  $\MDC(\Z\cap [5,20])$ indicating the heights of its elements.



\begin{lem} \label{MDC-lemma1} 
  Let $U$ be an interval in $\Z$.  Say $\MDC(U) = \{U_1, \ldots,
  U_n\}$ where the $U_i$ are indexed in ascending order.  Then for each $k$ either $$\left| \bigcup_{i=1}^{k-1} U_i \right|
\  < \  |U_k|  \hspace{1cm} \mathit{or} \hspace{1cm} \left| \bigcup_{i=k+1}^{n} U_i \right|
 \ < \ |U_k|.$$
\end{lem}

\bp Let  $r := \hght(U_k)$.  Assume $k$ is not $1$ or $n$, as otherwise the result is
trivial.  At most two
consecutive $U_i$ have the same height, so either $U_{k-1}$ or
$U_{k+1}$ is at a different height to $U_k$.  By the maximality of
the elements of $\MDC(U)$, if one of $U_{k-1}$ or $U_{k+1}$ has
height greater than $r$ then the other must have height less than $r$.  So $\hght(U_{k-1}) < r$ or $\hght(U_{k+1}) < r$.

As suggested by the example of Figure~\ref{fig:MDC}, $\MDC(U)$
resembles a pyramid, in that there is some $m$ such that the
sequence $\hght(U_1), \ldots, \hght(U_m)$ is strictly increasing and
the sequence $\hght(U_{m+1}), \ldots, \hght(U_n)$ is strictly
decreasing.  So if $\hght(U_{k-1}) < r$ then the sequence
$\hght(U_1), \ldots, \hght(U_{k-1})$ is strictly increasing and 
$$\left| \bigcup_{i=1}^{k-1} U_i \right| \ \leq  \ \sum_{i=0}^{r-1} 2^i \ <
\ 2^r  \ = \  |U_k|.$$  The case where $\hght(U_{k+1}) < r$ is similar. \ep

\bs

In order to control the merging process mentioned in the introduction,
we will use the following lemmas.  The first is an observation
on how dyadic covers can be converted to minimal dyadic covers.

\begin{lem}\label{lem:dyadicToMDC}
  Let $U \subseteq \Z$ be an interval and $X_0$ be a dyadic cover of
$U$.  Then there is a
  sequence $X_0, \ldots, X_n$ of
  dyadic covers of $U$ with  $X_n = \MDC(U)$ and, for all $i$, $X_{i+1}$ is obtained from $X_i$ by merging two adjacent dyadic intervals  --- that is,  if $X_i = \{S_1,
  \ldots, S_r\}$, where the $S_j$ are indexed in ascending order,
  then there exists $k$ such that $|S_k| = |S_{k+1}|$ and $X_{i+1} = \{S_1,
  \ldots, S_k \cup S_{k+1}, \ldots, S_r\}$.  In particular, if $S \in X_i$, then there exists $W \in X_0$ such that $W \subseteq S$.
\end{lem}

\bp If $X_0=\MDC(U)$, we are done.  Otherwise, by the definition of
$\MDC(U)$, there are dyadic intervals $W \in X_0$ and $D \subseteq U$
such that $W\subsetneq D$.  Choose a minimal length such $W \in X_0$.  The dyadic intervals containing
$W$ are well-ordered by containment; let $D$ be the
dyadic interval of minimal length strictly containing $W$.  We then have
$W \subsetneq D \subseteq U$
and $\ell(D)=2\ell(W)$.

We claim that $W' := D \smallsetminus W$  is in $X_0$.  We know that $W'$ is a dyadic
interval of the same length as $W$. Since $W' \subseteq U$ and
$X_0$ covers $U$, there must be an interval $Z \in X_0$ such that $Z \cap
W'$ is nonempty.  We claim that $Z=W'$.  Since $Z$ and
$W'$ are dyadic and have nonempty intersection, one must contain
the other.  This containment, however, cannot be strict; on one hand,
if $W' \subsetneq Z$, it must also contain $W$, so since $X_0$
is a collection of disjoint sets, $Z$ cannot strictly contain
$W'$.  On the other hand, by the minimality of $W$, the set
$Z$ cannot strictly contain $W'$.  Thus $W'=Z \in X_0$.

Then $W$ and $W'$ are adjacent intervals of equal length; without loss of generality, assume that $W$ is to the left of $W'$, so that $X_0$ can be expressed as
$$X_0 \  =  \ \{S_1,  \ldots, W, W', \ldots, S_r\}$$
Then we let 
$$X_1 \ = \ \{S_1,  \ldots, W \cup W', \ldots, S_r\}.$$
We repeat the process to construct $X_2,\dots,X_n$.  With each step, the number of elements in the partition decreases by one, and so the process terminates, and then every element of $X_n$ is maximal among dyadic intervals contained in $U$, so $X_n=\MDC(U)$.

The last assertion in the lemma follows by induction on $n$: if $n=0$,
it is trivially true, and by construction, any element of $X_{i+1}$
contains an element of $X_i$.  \ep

\bs

We can now describe a process of merging two adjacent
minimal dyadic covers.  Note that all of the changes occur at
the boundary between the two covers.

\begin{cor} \label{cor:mergingcovers1}
  Let $U, V \subseteq \Z$ be adjacent intervals with $\min \left(V\right) =
  \max \left( U\right) +1$.  Then there exists a
  sequence $X_0, \ldots, X_n$ of
  dyadic covers of $U \cup V$ with $$X_0 \ = \  \MDC(U) \cup \MDC(V)  \hspace{5mm} \text{and} \hspace{5mm} X_n \ =  \ \MDC(U \cup V)$$ and, for each $i$, expressing $X_i$ as $\{S_1,
  \ldots, S_r\}$  where the $S_j$ are indexed in ascending order,
   there exists $k$ such that $|S_k| = |S_{k+1}|$, $X_{i+1} = \{S_1,
  \ldots, S_k \cup S_{k+1}, \ldots, S_r\}$ and $S_k\cup S_{k+1}$ is not a subset of $U$ or $V$.
\end{cor}
\bp $\MDC(U) \cup \MDC(V)$ is a dyadic cover, so we may apply Lemma~\ref{lem:dyadicToMDC} to obtain a sequence $X_0,\dots, X_n=\MDC(U\cup
V)$ of dyadic covers.  It remains only to prove the final assertion.
By contradiction, suppose that $S_k\cup S_{k+1}\subset U$.  Then there
is some $W\in \MDC(U) \cup \MDC(V)$ such that $W\subseteq
S_k$.  In fact, we must have $W\in \MDC(U)$.  But  this is
impossible since $S_k\cup S_{k+1}$ is a dyadic interval strictly
containing $W$ and contained in $U$.  Similarly, $S_k\cup S_{k+1}$ cannot be a subset of $V$.\ep

\bs

The following similar observation applies when one of the intervals
consists of a single integer.

\begin{cor} \label{cor:mergingcovers2}
  Suppose $U \subset \Z$ is an interval.  Let $s = \min\left(U\right) - 1$ and
  $t = \max\left(U\right) + 1$.  Then: \begin{enumerate}
    \item There exists a sequence $X_0, \ldots, X_n$ of dyadic
    covers of $\{s\} \cup U$ with $$X_0 \ = \ \{ \{s\}\} \cup \MDC(U) \hspace{5mm}
    \mathit{and}  \hspace{5mm} X_n \ = \ \MDC(\{s\} \cup U)$$ and, for each $i$, expressing $X_i$ as $\{S_1, \ldots, S_r\}$  where the $S_i$ are indexed in ascending
    order,  we find $|S_1| = |S_2|$ and $X_{i+1} = \{S_1 \cup S_2, S_3, \ldots,
    S_r\}$.
    \item There exists a sequence $Y_0, \ldots, Y_n$ of dyadic
    covers of $U \cup \{t\}$ with $$Y_0 \ = \ \MDC(U) \cup \{\{t\} \}  \hspace{5mm}
    \mathit{and} \hspace{5mm} Y_n \ = \ \MDC(U \cup \{t\})$$ and, for each $i$, expressing $Y_i$ as $\{S_1, \ldots, S_r\}$  where the $S_i$ are indexed in ascending
    order, we find $|S_{r-1}| = |S_r|$ and $Y_{i+1} = \{S_1, \ldots,
    S_{r-2}, S_{r-1} \cup S_r\}$.
  \end{enumerate}
\end{cor}

\ss

\begin{defn}[Partitioned Alternating Form]\label{def:paf}
Let $v=v(a,b,c,d)$ be a balanced word partitioned as $v=v_1\ldots v_k$.
Let $\lambda_i$ \textup{(}resp.\ $\mu_i$\textup{)} be $v_i$ with all letters $c^{\pm 1}$ and $d^{\pm 1}$ \textup{(}resp.\   $a^{\pm 1}$ and $b^{\pm 1}$\textup{)} deleted.
Obtain $\rho_i(a,b,c)$ from $\lambda_i$ by replacing each $a, b, a^{-1},b^{-1}$  by $ac^{-1}, bc^{-1}, ca^{-1},cb^{-1}$, respectively.
Obtain $\sigma_i(a,c,d)$ from $\mu_i$ by replacing each $c,d, c^{-1},d^{-1}$ by  $ca^{-1}, da^{-1}, ac^{-1},ad^{-1}$, respectively.
The {\em partitioned alternating form} of $v$ with respect to  $v_1 \ldots v_k$ is
\[\tau \ := \ \rho_1\, (ca^{-1})^{\xi(\lambda_1)}\sigma_1\, (ac^{-1})^{\xi(\lambda_1\mu_1)}
\ldots
\rho_k \, (ca^{-1})^{\xi(\lambda_1\mu_1\ldots\lambda_k)}\sigma_k.\]
\end{defn}

\begin{lem}\label{partitioned lemma}
Any partitioned alternating form for $v$ equals $v$ in $F(a,b) \times F(c,d)$.
\end{lem}

\bp
We will use the notation of Definition~\ref{def:paf}.
Each $v_i= \lambda_i \mu_i$ in $S$ and so  $v= \lambda_1 \mu_1 \ldots \lambda_k \mu_k$ in $F(a,b) \times F(c,d)$.

Commutator relations can be used to convert a word $c^m \lambda \mu$, where $\lambda=\lambda(a,b)$, $\mu=\mu(c,d)$, and $m \in \mathbb{Z}$,  to $\rho \, (ca^{-1})^{m+ \xi(\lambda)} \,  \sigma \, (ac^{-1})^{m+ \xi \, (\lambda \mu)} \, c^{m+ \xi(\lambda \mu)}$,  where $\rho$ and  $\sigma$ are obtained from $\lambda$ and $\mu$ as per  Definition~\ref{def:paf}.
Making $k$ successive such transformations, working from left to right and beginning with $m=0$, converts $\lambda_1 \mu_1 \ldots \lambda_k \mu_k$ to $\tau$.  (Note that $\xi(\lambda_1\mu_1\ldots\lambda_k\mu_k)=0$ as $v$ is balanced.)
\ep

\begin{lem}\label{partitioned dependence}
In Definition~\textup{\ref{def:paf}}, for all $i$, the prefix
$$ \rho_1\, (ca^{-1})^{\xi(\lambda_1)}\sigma_1\, (ac^{-1})^{\xi(\lambda_1\mu_1)}
\ldots   \rho_{i} \, (ca^{-1})^{\xi(\lambda_1\mu_1\ldots\lambda_{i})} \, \sigma_{i} \,  (ac^{-1})^{\xi(\lambda_1\mu_1\ldots \lambda_{i}\mu_{i})}$$
and suffix
$$ \rho_{i+1}\, (ca^{-1})^{\xi(\lambda_1\mu_1 \ldots \lambda_{i+1})}\sigma_{i+1} \, (ac^{-1})^{\xi(\lambda_1\mu_1 \ldots \lambda_{i+1} \mu_{i+1} )}
\ldots   \rho_{k} \, (ca^{-1})^{\xi(\lambda_1\mu_1\ldots\lambda_{k})} \, \sigma_k  $$
of $\tau$ depend only on $v_1 \ldots v_i$ and $v_{i+1} \ldots v_k$, respectively.
\end{lem}

\bp
In the case of the prefix, this is self-evident.  It is true for the suffix, because  $\xi(\lambda_1\mu_1\ldots\lambda_k \mu_k)=0$ as $v$ is balanced, which allows one to express the exponents of $(ca^{-1})$ and $(ac^{-1})$ in terms of $\mu_{i+1}$,  $\lambda_{i+2}$, $\mu_{i+2}$, \ldots, $\lambda_k$, and $\mu_k$.
\ep

 \ms

\begin{defn}[Dyadic Alternating Form]\label{def:daf}
Let $v=v(a,b,c,d,s)$ be a balanced subword of a word $w$. Let
$\hat{v}$ and $\hat{w}$ be $v$ and $w$ with all $s^{\pm 1}$ removed.
Let $V$ be the set of positions of the letters of
$\hat{v}$ in $\hat{w}$.  Define the \emph{dyadic partition} of
$\hat{v}$ to be the partition corresponding to the partition $\MDC(V)$ of $V$.  The
\emph{dyadic alternating form} of $v$ is defined to be the
partitioned alternating form of $\hat{v}$ with respect to the dyadic
partition.
\end{defn}

Note that the dyadic alternating form of $v$ depends on its position in $w$.

\begin{lem} \label{dyadic equals}
The dyadic alternating form produced in the above definition equals $v$ in $S$.
\end{lem}

\bp
Since $v$ is balanced in this definition, $v= \hat{v}$ in $S$.  The result then follows from Lemma~\ref{partitioned lemma}.
\ep

\bs

One might think that the $(ac^{-1})^{\xi(\rho_1\ldots)}$ and $(ca^{-1})^{\xi(\rho_1\ldots)}$ inserted  could dramatically increase length, but the following estimates show this is not so for dyadic alternating form.

\begin{lem}\label{lem:exp-sum}
In the dyadic partition $v_1 \ldots v_k$ of $\hat{v}$ arising in Definition~\textup{\ref{def:daf}}, for all $i$,
$$|\xi(v_1\ldots v_i)| <2\ell(v_i) \hspace{5mm} \text{and}  \hspace{5mm} |\xi(v_1\ldots v_{i-1}\lambda_i)| <2\ell(v_i).$$
\end{lem}

\bp As the $v_i$ are defined using a minimal dyadic cover, Lemma
\ref{MDC-lemma1} implies that either $\ell(v_1\ldots
v_{i-1}) <\ell(v_i)$ or $\ell(v_{i+1}\ldots v_k)<\ell(v_i)$.

In the first case,  $|\xi(v_1\ldots v_i)| \leq \ell(v_1\ldots
v_i)=\ell(v_1\ldots v_{i-1})+\ell(v_i)<2\ell(v_i)$ and
$|\xi(v_1\ldots v_{i-1}\lambda_i)| \leq \ell(v_1\ldots
v_{i-1}\lambda_i)=\ell(v_1\ldots
v_{i-1})+\ell(\lambda_i)<2\ell(v_i)$.

In the second case,  $\xi(v_{i+1}\ldots v_k)=-\xi(v_1\ldots v_i)$ since $\xi(\hat{v})=\xi(v)=0$ as $v$ is balanced.
So  $|\xi(v_1\ldots v_i)|=|\xi(v_{i+1}\ldots v_k)|\leq \ell(v_{i+1}\ldots v_k)<\ell(v_i)$. Similarly, $\xi(\mu_iv_{i+1}\ldots v_k)=-\xi(v_1\ldots v_{i-1}\lambda_i)$ and
so  $|\xi(v_1\ldots v_{i-1}\lambda_i)|=|\xi(\mu_iv_{i+1}\ldots v_k)\leq \ell(\mu_iv_{i+1}\ldots v_k)<2\ell(v_i)$.
\ep

\begin{lem}\label{lem:mdaf-ten-n}
The dyadic alternating form of $v$ defined in Definition~\textup{\ref{def:daf}} has  length at most $10\ell(v)$.
\end{lem}

\bp In the dyadic alternating form, each letter of $\hat{v}$ is matched with either an $a^{\pm 1}$ or $c^{\pm 1}$.  As $\ell(\hat{v})\leq \ell(v)$, this accounts for at most $2\ell(v)$ letters.

There are $2(|\xi(\lambda_1)|+|\xi(\lambda_1\mu_1)|+\ldots+|\xi(\lambda_1\mu_1\ldots\lambda_k)|)$ further letters appearing in the powers of $(ac^{-1})^{\pm 1}$.
But for all $i$,
\begin{equation*}
\begin{split}
|\xi(\lambda_1\mu_1&\ldots \lambda_{i-1}\mu_{i-1}\lambda_i)|+|\xi(\lambda_1\mu_1\ldots \lambda_{i}\mu_{i})|  \\ &= \ |\xi(v_1\ldots v_{i-1}\lambda_i)|+|\xi(v_1\ldots v_i)|  \ < \  4\ell(v_i)
\end{split}
\end{equation*}
 by  Lemma~\ref{lem:exp-sum}.  So it suffices to add
 $$2(|\xi(\lambda_1)|+|\xi(\lambda_1\mu_1)|+\ldots+|\xi(\lambda_1\mu_1\ldots\lambda_k)|) \  \leq \
 \sum_{i=1}^k 8\ell(v_i) \  = \  8\ell(\hat{v}) \ \leq \  8\ell(v).$$
\vspace*{-5mm}
\ep

\section{Subroutines} \label{sec:subroutines}

In this section we present the key subroutines that will be called by
the main algorithm.  Algorithms 1--6 manipulate words using the
relations of $F(a, b) \times F(c, d)$, while Algorithm \ref{alg:main}
uses the relators of the whole group $S$.  Each algorithm gives a
method for converting one word to another word by applications of
relators, free expansions, and free reductions.  Thus each algorithm
has a cost in the sense of Definition \ref{cost etc}.

The first algorithm gives a rough-and-ready scheme for converting
words when one does not have any additional information about the
structure of the input word.  As such the cost of the algorithm is
high.

\ms

  \begin{algorithmbis}[Shuffling between words  $F(a,b)\times F(c,d)$]\label{alg:FxF}

\bs

\parbox{119.5mm}{\hspace*{-3mm}\parbox{17mm}{\textbf{Input:}}    Words $u=u(a,b,c,d)$ and $v=v(a,b,c,d)$ representing the same
elements of $F(a,b)\times F(c,d)$

\bs

\hspace*{-3mm}\parbox{17mm}{\textbf{Goal:}} Convert $u$ to $v$.

\bs

\hspace*{-3mm}\parbox{17mm}{\textbf{Method:}} Shuffle the $a$'s and $b$'s in $u$ to the front of the word and
freely reduce to produce a word $w$.
 The same procedure would convert  $v$ to $w$, so run it in reverse to
convert $w$ to $v$.}

\vspace*{8mm}

   \end{algorithmbis}

\ms
 
\begin{lem}\label{lem:FxF}
The cost \textup{(}in the sense of Definition~\textup{\ref{cost etc}}\textup{)} of the transformation of Algorithm~\textup{\ref{alg:FxF}} is at most $\ell(u)^2+ \ell(v)^2$.
\end{lem}

\bp
Each letter in $u$ or $v$ is shuffled past fewer than $\ell(u)$ or $\ell(v)$ other letters (respectively).
\ep

\bs  

  \begin{algorithmbis}[Merging within partitioned alternating form]\label{alg:small merge}

\bs

\parbox{119.5mm}{\hspace*{-3mm}\parbox{17mm}{\textbf{Input:}}   A balanced word $v = v(a,b,c,d)$, a partition $v_1 \ldots v_k$  of $v$, an  integer   $ i \in \{ 1, \ldots, k-1\}$, and the partitioned  alternating form  $\tau$ of $v$ with respect to $v_1 \ldots v_k$}

\bs

\parbox{119.5mm}{\hspace*{-3mm}\parbox{17mm}{\textbf{Goal:}}  Convert $\tau$ to the  
  partitioned alternating form  $\bar{\tau}$ of $v$ with respect to the partition
  $v_1 \ldots v_{i-1} \, \bar{v} \,  v_{i+2} \ldots v_k$
  where $\bar{v} = v_i v_{i+1}$.} 
 
\bs

\parbox{119.5mm}{\hspace*{-3mm}\parbox{17mm}{\textbf{Method:}} 
  We use the notation of Definition~\ref{def:paf} and
  write $$\tau \  = \ \rho_1\,
  (ca^{-1})^{\xi(\lambda_1)}\sigma_1\, (ac^{-1})^{\xi(\lambda_1\mu_1)}
  \ldots \rho_k \,
  (ca^{-1})^{\xi(\lambda_1\mu_1\ldots\lambda_k)}\sigma_k.$$
  Express $\tau$ as $\tau^{(0)} \tau^{(1)} \tau^{(2)}$ where
  $$\hspace*{-4mm}\begin{array}{lll}
   \tau^{(1)} & =   & \rho_i \, (ca^{-1})^{\xi(\lambda_1\mu_1\ldots \lambda_i)}\sigma_i
   \, (ac^{-1})^{\xi(\lambda_1\mu_1 \ldots \lambda_i\mu_i)}   \\
  & & \ \ \   \rho_{i+1} \, (ca^{-1})^{\xi(\lambda_1\mu_1\ldots\lambda_{i+1})} \,
  \sigma_{i+1} \, (ac^{-1})^{\xi(\lambda_1\mu_1\ldots\lambda_{i+1} \mu_{i+1})}.
  \end{array}$$
  Let $$\bar{\tau}^{(1)} \ = \ \bar{\rho} \, (ca^{-1})^{\xi(\lambda_1\mu_1\ldots
  \lambda_{i-1} \mu_{i-1} \bar{\lambda} )} \bar{\sigma} \,
  (ac^{-1})^{\xi(\lambda_1\mu_1 \ldots \lambda_{i-1}\mu_{i-1}
  \bar{\lambda} \bar{\mu} )}$$
  where $\bar{ \lambda}$ (resp.\ $\bar{\mu}$) is $\bar{v}$
  with all $c^{\pm 1}$ and $d^{\pm 1}$ (resp.\ $a^{\pm 1}$ and
  $b^{\pm 1}$) removed.  (So $\bar{ \lambda} = \lambda_{i}\lambda_{i+1}$
  and $\bar{ \mu} = \mu_{i}\mu_{i+1}$.)
  It follows from Lemma~\ref{partitioned dependence} that
  $\bar{\tau} =  \tau^{(0)} \bar{\tau}^{(1)} \tau^{(2)}$.  So  $\tau^{(1)}
   = \bar{\tau}^{(1)}$ in $F(a,b) \times F(c,d)$.} 
 
 \ms
 
\parbox{119.5mm}{Obtain $\bar{\tau}$ from $\tau$ by changing the subword
  $\tau^{(1)}$ to $\bar{\tau}^{(1)}$ using Algorithm~\ref{alg:FxF}.}

\vspace*{6mm}

     \end{algorithmbis}

\ms

\begin{lem}\label{small merge cost}
The cost of Algorithm~\textup{\ref{alg:small merge}} is at most $136 \left( \ell(v_i v_{i+1}) + | \xi(v_1 \ldots v_{i-1} )| \right)^2$.
\end{lem}

\bp
The length of $\tau^{(1)}$ is at most
\begin{equation*}
\begin{split}
\ell(\rho_i \sigma_i \rho_{i+1} \sigma_{i+1})
\ +  \ & 2 \left( \rule{0mm}{4mm}  | \xi( \lambda_1 \mu_1 \ldots \lambda_i) |
\ + \  | \xi( \lambda_1 \mu_1 \ldots \lambda_i \mu_i) | \right. \\ & \ \ \ \ \  \left.
\ + \  | \xi( \lambda_1 \mu_1 \ldots \lambda_{i+1}) |
\ + \  | \xi( \lambda_1 \mu_1 \ldots \lambda_{i+1} \mu_{i+1}) |  \rule{0mm}{5mm} \right).
\end{split}
\end{equation*}
As $\ell(\rho_i \sigma_i \rho_{i+1} \sigma_{i+1}) = 2\ell(\lambda_i \mu_i \lambda_{i+1} \mu_{i+1}) = 2\ell(v_i v_{i+1})$ and each of the four other terms differs from $| \xi(\lambda_1 \mu_1 \ldots \lambda_{i-1} \mu_{i-1} )| = | \xi(v_1 \ldots v_{i-1} )|$ by at most $\ell(\lambda_i \mu_i \lambda_{i+1} \mu_{i+1})$,
$$\ell(\tau^{(1)}) \ \leq \  10\ell(v_i v_{i+1}) +  8 | \xi(v_1 \ldots v_{i-1} )|.$$
Similarly, $$\ell(\bar{\tau}^{(1)}) \ \leq \ 6\ell(v_i v_{i+1}) +  4
| \xi(v_1 \ldots v_{i-1} )|.$$ By Lemma~\ref{lem:FxF} the cost of
Algorithm~\ref{alg:small merge} is at most
$\ell(\tau^{(1)})^2+\ell(\bar{\tau}^{(1)})^2$ ---
this then gives the (crude) estimate we claim. \ep

\bs

Our next subroutine merges two subwords in dyadic alternating form into
one.

\newpage

  \begin{algorithmbis}[Merge two dyadic alternating subwords]\label{alg:bigmerge}

\bs

\parbox{119.5mm}{\hspace*{-3mm}\parbox{17mm}{\textbf{Input:}}   Two balanced words $u$ and $v$, such that $uv$ is a subword of 
$w = w(a,b,c,d,s)$, with $\tau_u$ and $\tau_v$
    the dyadic alternating forms of $u$ and $v$ with respect to their positions in $w$. 
  
\bs

\hspace*{-3mm}\parbox{17mm}{\textbf{Goal:}}   Convert $\tau_u \tau_v$ to the dyadic    alternating form $\tau$ of $uv$.

\bs

\hspace*{-3mm}\parbox{17mm}{\textbf{Method:}} Note that  $uv$ is balanced as it is the
    concatenation of two balanced words.   
  Let $\hat{u}$,  $\hat{v}$ and $\hat{w}$ be $u$, $v$ and $w$
    respectively with all occurrences of the letters $s^{\pm1}$
    removed.  Let $U, V \subset \Z$ be the sets of positions of
    the letters of $\hat{u}$ and $\hat{v}$ respectively in
    $\hat{w}$.  Let $X_0, \ldots, X_n$ be a sequence of dyadic
    covers of $U \cup V$ as given by Corollary~\ref{cor:mergingcovers1}.
    Let $\sigma_i$ be the partition of $\hat{u}\hat{v}$
    corresponding to $X_i$ and let $\tau_i$ be the partitioned
    alternating form of $\hat{u}\hat{v}$ with respect to $\sigma_i$.
    Then $\tau_0 = \tau_u \tau_v$ and $\tau_n = \tau$.  For each
    $i$, convert $\tau_i$ to $\tau_{i+1}$ by calling Algorithm~\ref{alg:small merge}.
}

\vspace*{10mm}

   \end{algorithmbis}

\bs

We postpone a full cost analysis of this algorithm to
Section~\ref{sec:cost}.  In fact, there we will estimate the total
cost of all the calls of Algorithm~\ref{alg:small merge} throughout
our main algorithm rather than their total cost within any single
call on  Algorithm~\ref{alg:bigmerge}.  However, we will pause to
give the following lemma which will be crucial to that analysis.

\begin{lem} \label{key cost lemma} 
We continue with the notation of Algorithm~\textup{\ref{alg:bigmerge}}.  Fix
$i$ and say that the partition $\sigma_{i+1}$ is formed from
$\sigma_i$ by combining the two subwords $\theta$ and $\phi$.  Then
$\ell(\theta) = \ell(\phi)$ and either $\theta$ is a subword of
$\hat{u}$ or $\phi$ is a subword of $\hat{v}$.  Furthermore the cost
of the call to Algorithm~\textup{\ref{alg:small merge}} converting $\tau_i$
to $\tau_{i+1}$ is at most $2176 \, \ell(\theta)^2$.
\end{lem}

\bp 
Say $X_i = \{S_1, \ldots, S_r\}$ and $X_{i+1} = \{S_1, \ldots, S_k
\cup S_{k+1}, \ldots, S_r\}$, where the indexing on the $S_i$ is
ascending.  Say that the partition $\sigma_i$ of $\hat{u}\hat{v}$ is
$\pi_1 \ldots \pi_r$. Then $\theta = \pi_k$ and $\phi = \pi_{k+1}$. By
Corollary~\ref{cor:mergingcovers1}, $\ell(\pi_k) = \ell(\pi_{k+1})$ and
$\pi_k\pi_{k+1}$ is not a subword of either $\hat{u}$ or $\hat{v}$.

By Lemma~\ref{small merge cost}, the cost of the call to
Algorithm~\ref{alg:small merge} in question is at most
$136(\ell(\pi_k\pi_{k+1}) + |\xi(\pi_1 \ldots \pi_{k-1})|)^2$.  Note
that $\pi_1 \ldots \pi_{k-1}$ is a subword of $\hat{u}$, which in turn
is a subword of $\pi_1 \ldots \pi_{k+1}$, and thus
$$|\xi(\pi_1 \ldots \pi_{k-1})-\xi(\hat{u})| \ \le \ \ell(\pi_k\pi_{k+1}) \ = \ 2\ell(\pi_k).$$
Furthermore, since $\hat{u}$ is balanced, $\xi(\hat{u})=0$.  Thus $|\xi(\pi_1 \ldots
\pi_{k-1})|\le 2\ell(\pi_k)$ and
$$136(\ell(\pi_k\pi_{k+1}) + |\xi(\pi_1 \ldots \pi_{k-1})|)^2 \ \leq \
136(2\ell(\pi_k) + 2\ell(\pi_k))^2 \ \leq \ 2176 \ell(\pi_k)^2.$$ \ep



\bs

  \begin{algorithmbis}[Shuffling and canceling a $c^{\pm 1}$]\label{alg:pass_c}

\bs

\parbox{119.5mm}{\hspace*{-3mm}\parbox{17mm}{\textbf{Input:}}     A word 
 $\tau$ of the form $$\rho_1\, (ca^{-1})^{\alpha_1}\sigma_1\,
    (ac^{-1})^{\beta_1} \ \ldots   \    \rho_{k} \, (ca^{-1})^{\alpha_k}\sigma_{k}
    (ac^{-1})^{\beta_k}$$ where $\alpha_i, \beta_i \in \mathbb{Z}$,
    $\rho_i=\rho_i(a,b,c)$ and $\sigma_i=\sigma_i(a,c,d)$.  Let
    $\epsilon \in \{\pm1\}$ and define $$\bar{\tau} \ = \ \rho_1\,
    (ca^{-1})^{\alpha_1+\epsilon}\sigma_1\, (ac^{-1})^{\beta_1+\epsilon} \
    \ldots   \    \rho_{k} \, (ca^{-1})^{\alpha_k+\epsilon}\sigma_{k}
    (ac^{-1})^{\beta_k+\epsilon}.$$

\bs

\hspace*{-3mm}\parbox{17mm}{\textbf{Goal:}} Convert $c^{\epsilon} \tau  c^{-\epsilon}$ to $\bar{\tau}$.

\bs

\hspace*{-3mm}\parbox{17mm}{\textbf{Method:}}  Working from left to right, shuffle $c^{\epsilon}$ through $\rho_1$,
    replace $c^{\epsilon}$ by $(ca^{-1})^{\epsilon} a^{\epsilon}$,
    shuffle  $a^{\epsilon}$ through $(ca^{-1})^{\alpha_1}\sigma_1$,
    replace $a^{\epsilon}$ by $(ac^{-1})^{\epsilon} c^{\epsilon}$,
    and continue similarly.  When  $c^{\epsilon}$ emerges after
    $(ac^{-1})^{\beta_k}$ cancel it with the $c^{-\epsilon}$.
}

\vspace*{10mm}

   \end{algorithmbis}

\bs

\begin{lem} \label{pass_c cost}
The cost of Algorithm~\textup{\ref{alg:pass_c}} is at most $2k + \ell(\tau)$.
\end{lem}

\bp
Replacing $c^{\epsilon}$ by $(ca^{-1})^{\epsilon} a^{\epsilon}$ costs $1$ when $\epsilon = -1$ and costs $0$ otherwise.   The same is true of replacing  $a^{\epsilon}$ by $(ac^{-1})^{\epsilon} c^{\epsilon}$.  The other contributions to cost stem from carrying   $a^{\epsilon}$ or $c^{\epsilon}$ past letters of $\tau$.
\ep

\bs

Our next subroutine expands a dyadic form subword by assimilating a
letter on each side to produce a word in partitioned (but not
necessarily dyadic) alternating form.

\bs

  \begin{algorithmbis}[Assimilate $x,y$ into a dyadic alternating word: I]\label{alg:assimilate1}

\bs

\parbox{119.5mm}{\hspace*{-3mm}\parbox{17mm}{\textbf{Input:}}    A subword of 
  $w = w(a, b, c, d, s)$  of
    the form $xvy$ where $x, y \in \{a, b, c, d\}^{\pm1}$ have
    opposite exponents and $v = v(a, b, c, d, s)$ is balanced with dyadic alternating form 
    $\tau$.
    Write $\hat{v}$ and $\hat{w}$ for the words $v$ and $w$
    respectively with all occurrences of the letter $s^{\pm1}$
    removed.  Say $\hat{v}$ has dyadic partition $v_1 \ldots v_k$ in
    $\hat{w}$.} 

\bs

\parbox{119.5mm}{\hspace*{-3mm}\parbox{17mm}{\textbf{Goal:}}     Convert $x \tau y$ into  partitioned alternating form $\bar{\tau}$
 with respect to the partition $xv_1 \ldots v_k y$.}

\bs

\parbox{119.5mm}{\hspace*{-3mm}\parbox{17mm}{\textbf{Method:}}     Using the notation of Definition~\ref{def:paf},
    $$\tau \ = \  \rho_1\, (ca^{-1})^{\xi(\lambda_1)}\sigma_1\, (ac^{-1})^{\xi(\lambda_1\mu_1)}
    \ldots \rho_k \, (ca^{-1})^{\xi(\lambda_1\mu_1\ldots\lambda_k)}\sigma_k.$$}

\parbox{119.5mm}{Write $x$ and $y$ as $\lambda_0 \mu_0$ and  $\lambda_{k+1}
    \mu_{k+1}$, respectively, where $\lambda_0$ and $\lambda_{k+1}$ are
    each in $\{ a, a^{-1}, b, b^{-1}, \e \}$ and  $\mu_0$ and
    $\mu_{k+1}$ are each in $\{ c, c^{-1}, d, d^{-1}, \e \}$.}

\parbox{119.5mm}{Thus
    $$\begin{array}{lll}
    \bar{\tau} & = &   \rho_0\, (ca^{-1})^{\xi(\lambda_0)}\sigma_0 \, (ac^{-1})^{\xi(\lambda_0\mu_0)}
    \ldots \rho_{k+1} \, (ca^{-1})^{\xi(\lambda_0\mu_0\ldots\lambda_{k+1})}\sigma_{k+1} \\
    & =\rule{0mm}{5mm} &  \rho_0\, (ca^{-1})^{\xi(\lambda_0)}\sigma_0 \, (ac^{-1})^{\xi(x)} \rho_1 \, (ca^{-1})^{\xi(x) + \xi(\lambda_1)}\sigma_1 \, (ac^{-1})^{\xi(x) + \xi(\lambda_1 \mu_1)} \\
     & & \ \ \ \ldots   \rho_{k+1} \, (ca^{-1})^{\xi(x) + \xi(\lambda_1\mu_1\ldots\lambda_{k+1})}\sigma_{k+1}.
    \end{array}$$}

\parbox{119.5mm}{Transform  $x \tau y$ to $\bar{\tau}$ by applying Algorithm~\ref{alg:FxF} to convert
    $$\begin{array}{lll}
     x & \mapsto & \rho_0 (ca^{-1})^{\xi(\lambda_0)} \sigma_0 (ac^{-1})^{\xi(\lambda_0 \mu_0)} c^{\xi(x)} \ \ \textup{and} \\
     y & \mapsto\rule{0mm}{5mm} & c^{\xi(y)} \rho_{k+1}(ac^{-1})^{\xi(\mu_{k+1})} \sigma_{k+1}, \\
        & =\rule{0mm}{5mm} & c^{-\xi(x)} \rho_{k+1}(ca^{-1})^{\xi(x) + \xi(\lambda_1\mu_1\ldots\lambda_{k+1})} \sigma_{k+1},
         \end{array} $$
    and then applying Algorithm~\ref{alg:pass_c} to the subword $c^{\xi(x)} \tau c^{-\xi(x)}$.
}

\vspace*{10mm}

   \end{algorithmbis}

\ms

 \begin{lem} \label{easy little lemma}
 In $F(a,b) \times F(c,d)$,
 $$\begin{array}{lll}
 x & = & \rho_0 (ca^{-1})^{\xi(\lambda_0)} \sigma_0 (ac^{-1})^{\xi(\lambda_0 \mu_0)} c^{\xi(x)} \ \  \textup{ and }  \\
 y & =\rule{0mm}{5mm} & c^{\xi(y)} \rho_{k+1}(ac^{-1})^{\xi(\mu_{k+1})} \sigma_{k+1},
 \end{array}$$
 in the notation of Algorithm~\textup{\ref{alg:assimilate1}}.
  Moreover, the cost of transforming $x$ or $y$ in this way is at most $2$.
 \end{lem}

\bp
This is easily checked case-by-case.  The reason the cost is so low is that most of the moves involved are free-expansions rather than applications of commutator relations.
\ep

\begin{lem} \label{lem:cost_alg_assimilate1} 
  We continue with the notation of Algorithm~\textup{\ref{alg:assimilate1}}.
  The total cost of calling this algorithm is at most $3\ell(\tau) +
  4$.
\end{lem}

\bp
  By Lemmas~\ref{pass_c cost} and \ref{easy little lemma} the
  cost is at most $2k + \ell(\tau) +4$.  But $k \leq \ell(\tau)$.\ep

\bs

The following routine builds on Algorithm~\ref{alg:assimilate1}.  It
assimilates a letter on either side of a dyadic subword to produce
the dyadic form of the concatenated subword.

\bs

  \begin{algorithmbis}[Assimilate $x,y$ into a dyadic alternating word: II]\label{alg:assimilate2}

\bs

\parbox{119.5mm}{\hspace*{-3mm}\parbox{17mm}{\textbf{Input:}}     The partitioned alternating form 
 $\bar{\tau}$ of $xvy$ from Algorithm~\ref{alg:assimilate1}.
}

\bs

\parbox{119.5mm}{\hspace*{-3mm}\parbox{17mm}{\textbf{Goal:}} Convert the partitioned alternating form 
 $\bar{\tau}$ of $xvy$ from the previous algorithm to dyadic alternating form $\tau'$.}

\bs

\parbox{119.5mm}{\hspace*{-3mm}\parbox{17mm}{\textbf{Method:}} First apply Algorithm~\ref{alg:assimilate1} to $\tau$.  Let
  $\bar{\tau}$ be the output.}

 \bs

\parbox{119.5mm}{Let $\hat{v}$ and $\hat{w}$ be the words $v$ and $w$ respectively
  with all occurrences of the letters $s^{\pm1}$ deleted and let $V$
  be the set of positions of the letters of $\hat{v}$ in $\hat{w}$.
  Let $s = \min\left(V\right) - 1$ and $t = \max\left(V\right) + 1$.  Thus $s$ and $t$ are the
  positions of the letters $x$ and $y$ respectively in $\hat{w}$.}

 \bs

\parbox{119.5mm}{Let $X_0, \ldots, X_n$ be a sequence of dyadic covers of $\{s\}
  \cup V$, as given by part~1 of Corollary~\ref{cor:mergingcovers2},
  converting $X_0 = \{\{s\}\} \cup \MDC(V)$ to $X_n = \MDC(\{s\} \cup V)$.  For
  each $i$, define $Y_i$ to be $X_i \cup \{\{t\}\}$.  Thus $Y_0, \ldots, Y_n$ is a
  sequence of dyadic covers of $\{s\} \cup V \cup \{t\}$ converting
  $Y_0 = \{\{s\}\} \cup \MDC(V) \cup \{\{t\}\}$ to $Y_n = \MDC(\{s\}
  \cup V) \cup \{\{t\} \}$.}

 \bs

\parbox{119.5mm}{Let $Y_n, \ldots, Y_m$ be a sequence of dyadic
  covers of $\{s\} \cup V \cup
  \{t\}$, as given by part~2 of Corollary~\ref{cor:mergingcovers2},
  converting, converting $Y_n = \MDC(\{s\} \cup V) \cup \{\{t\}\}$
  to $Y_m = \MDC(\{s\} \cup V \cup \{t\})$.}

 \bs

\parbox{119.5mm}{For each $i = 0, \ldots, m$, let $\sigma_i$ be the partition of $x \hat{v}y$
  corresponding to $Y_i$, and let $\tau_i$ be the partitioned
  alternating form of $x \hat{v} y$ with respect to $\sigma_i$.
  Then $\tau_0 = \bar{\tau}$ and $\tau_m = \tau'$.  For each $i$,
  call Algorithm~\ref{alg:small merge} to convert $\tau_i$ to
  $\tau_{i+1}$.
}

\vspace*{10mm}

   \end{algorithmbis}

\bs

A lemma analogous to Lemma~\ref{key cost lemma} will be important when we come to analyse cost.

\begin{lem} \label{other key cost lemma} 
  We continue with the notation of Algorithm~\textup{\ref{alg:assimilate2}}.
  Fix $i$ and say the partition $\sigma_{i+1}$ of $x \hat{v} y$ is
  formed from $\sigma_i$ by combining the subwords $\theta$ and
  $\phi$.  Then $\ell(\theta) = \ell(\phi)$ and the cost of the call
  to Algorithm~\textup{\ref{alg:small merge}} which converts $\tau_i$ to
  $\tau_{i+1}$ is at most $2176 \ell(\theta)^2$.
\end{lem}

\bp
  First consider the case that $0 \leq i \leq n-1$.  Say that the
  partition $\sigma_i$ of $x \hat{v} y$ is $\pi_1 \ldots \pi_k$.
  Then $\theta = \pi_1$ and $\phi = \pi_2$ and the partition
  $\sigma_{i+1}$ is $\bar{\pi} \pi_3 \ldots \pi_k$, where
  $\bar{\pi} = \theta \phi$.  By part~1 of
  Corollary~\ref{cor:mergingcovers2}, $\ell(\theta) = \ell(\phi)$ and by
  Lemma~\ref{small merge cost} the cost of the call to Algorithm~\ref{alg:small
  merge} is at most $$136\ell(\theta\phi)^2 \ = \  136(2\ell(\theta))^2 \  = \ 
  544 \ell(\theta)^2.$$

  Now consider the case that $n \leq i \leq m-1$.  Say that the partition
  $\sigma_i$ of $x \hat{v} y$ is $\pi_1 \ldots \pi_k$.  Then $\theta
  = \pi_{k-1}$ and $\phi = \pi_k$ and the partition $\sigma_{i+1}$
  is $\pi_1 \ldots \pi_{k-1} \bar{\pi}$ where $\bar{\pi} = \pi_{k-1}
  \pi_k$.  By part~2 of Corollary~\ref{cor:mergingcovers2}, $\ell(\theta) = \ell(\phi)$ and by
  Lemma~\ref{small merge cost} the cost of the call to Algorithm~\ref{alg:small
  merge} is at most \begin{align*}
    136(\ell(\theta\phi) + |\xi(\pi_1 \ldots
  \pi_{k-1})|)^2 & \  =  \ 136(2\ell(\theta) + |\xi(\theta\phi)|)^2 \\
  & \ \leq \ 136(2\ell(\theta) + \ell(\theta\phi))^2 \\
  & \ = \ 136(2\ell(\theta) + 2\ell(\theta))^2 \\
  & \ = \ 2176 \ell(\theta)^2.\end{align*}\vspace*{-3mm}
\ep

   \section{Our main algorithm}
   \label{sec:main}

We are now ready to give our main algorithm, which converts a
null-homotopic word of length $n$  in $S$ to $\e$ via $n/2$
intermediate words, each of length at most $10n$, by applying
relations from the presentation~(\ref{pres}).

\ms

  \begin{algorithmbis}[Our main algorithm]\label{alg:main}

\bs

\parbox{119.5mm}{\hspace*{-3mm}\parbox{17mm}{\textbf{Input:}}  A word $w = w(a,b,c,d,s)$  representing $1$ in $S$

\bs

\hspace*{-3mm}\parbox{17mm}{\textbf{Goal:}} Reduce $w$ to $\e$ by applying defining relators of $S$.

\bs

\hspace*{-3mm}\parbox{17mm}{\textbf{Method:}} Define  $n := \ell(w)$.
As all the defining relators in the presentation \eqref{pres}  of $S$ are of even length, $n$ is even.
We will obtain a sequence of words $w_i$ and a sequence $T_i$ of subsets of $\{1, \ldots, n \}$ beginning with $w_0:=w$ and $T_0 := \emptyset$.}

\bs

\parbox{119.5mm}{In fact, $w$ and $T_i$ will define $w_i$ as follows.  Express $w$ as
\begin{eqnarray} \label{w} 
w & = & u_1v_1u_2\ldots v_{k-1}u_k,
\end{eqnarray}
where $T_i$ is the set of positions of the letters of the $v_j$ in $w$ and $v_1,u_2, \ldots, $ $u_{k-1}, v_{k-1} \neq \e$.  One sees inductively from the construction of successive $T_i$ below that the $v_j$ are balanced.
Let $\tau_j$ be the dyadic alternating form of  $v_j$ in $w$.  Then
\begin{eqnarray} \label{w_i}
w_i & = & u_1\tau_1u_2 \ldots \tau_{k-1} u_k.
\end{eqnarray}
}

\bs

\parbox{119.5mm}{We will now explain how to obtain $T_{i+1}$ from $T_i$ and then how to use relations in $S$ to transform $w_i$, expressed as (\ref{w_i}), to $w_{i+1}$.}

\bs

\parbox{119.5mm}{As $w$ and the $v_j$  are all balanced,  $u_1u_2\ldots u_k$ is also balanced by Lemma~\ref{lem:balanced-alternating}.  So Lemma~\ref{lem:xy} applies and tells us there is a subword $xy$ in $u_1\ldots u_n$ that either equals $s^{\pm 1}s^{\mp 1}$ or is such that $x$ and $y$ are in $\{a, b, c,  d \}^{\pm 1}$ and have opposite exponents.  Add the positions of $x$ and $y$ in $w$ to  $T_i$ to obtain $T_{i+1}$.}

\bs

\parbox{119.5mm}{In $w_i$ these $x$ and $y$ are either adjacent or separated by some $\tau_j$.}

\bs

\parbox{119.5mm}{If $xy=s^{\pm 1}s^{\mp 1}$, then remove $x$ and $y$  by shuffling $x$ through $\tau_j$ (if present) and cancelling it with $y$.}

\bs

\parbox{119.5mm}{If $x,y\in \{a,  b,  c,  d  \}^{\pm 1}$, then apply
Algorithm~\ref{alg:assimilate2} to replace $x\tau_jy$ with the
dyadic alternating form $\tau'$ of the subword $x v_jy$ of $w$.}

\bs

\parbox{119.5mm}{Next, if (in either case) we have brought two or three dyadic alternating form subwords together (that is, either $u_j=x$ and $j\neq 1$, or $u_{j+1}=y$ and $j+1 \neq k$, or both) then merge them using Algorithm~\ref{alg:bigmerge} (once or twice, as necessary).  The result is the word $w_{i+1}$.}

\bs

\parbox{119.5mm}{After $n/2$ iterations, $T_{n/2} = \{ 1, \ldots , n \}$ and $w_{n/2}$ contains no $s^{\pm 1}$   and represents $1$ in $F(a,b) \times F(c,d)$. Reduce $w_{n/2}$ to $\e$ using Algorithm~\ref{alg:FxF}.
}

\vspace*{8mm}

   \end{algorithmbis}

\section{Cost analysis}\label{sec:cost}

We will estimate the cost of our main algorithm in terms of $n=
\ell(w)$. Cost is incurred in four ways: shuffling an $s^{\pm 1}$ to
be cancelled with an $s^{\mp 1}$,   calls of
Algorithm~\ref{alg:bigmerge},  calls of
Algorithm~\ref{alg:assimilate2}, and  converting $w_{n/2}$ to $\e$
at the end.

Shuffling an $s^{\pm 1}$ to be cancelled with an $s^{\mp 1}$ costs at most $5n$ because the $\tau_j$ between them (if present) is alternating and so  $\ell(\tau_j)/2$ relations are required, and $\ell(\tau_j)/2 \leq 5\ell(v_j)\leq 5n$ by Lemma~\ref{lem:mdaf-ten-n}.
Since we do such shuffling at most $n/2$ times, this contributes no more than $5n^2/2$ to the total cost.

The final step,  converting $w_{n/2}$ to $\e$, costs at most $\ell(w_{n/2})^2  \leq 100 n^2$ by Lemmas~\ref{lem:mdaf-ten-n} and \ref{lem:FxF}.

Consider the call to Algorithm~\ref{alg:assimilate2}
converting $x \tau_j y$ to the dyadic alternating form of $x v_j y$.
The cost of this call can be divided into the cost of calls to
Algorithm~\ref{alg:assimilate1} and the cost of calls to
Algorithm~\ref{alg:small merge}.  By
Lemma~\ref{lem:cost_alg_assimilate1} each call to
Algorithm~\ref{alg:assimilate1} costs at most $3\ell(\tau_j) +4$ and
by Lemma~\ref{lem:mdaf-ten-n} we have that $\ell(\tau_j) \leq 10n$.
Thus the cost of each call to Algorithm~\ref{alg:assimilate1} is at
most $30n + 4$. In total Algorithm~\ref{alg:assimilate2}, and hence
Algorithm~\ref{alg:assimilate1}, is called at most $n/2$ times and
so the total cost of all calls to Algorithm~\ref{alg:assimilate1} is
at most $15n^2 + 2n$.

This leaves just the calls to Algorithm~\ref{alg:small merge} to
consider.  We will bound the total cost of all calls to this
algorithm arising from either Algorithm~\ref{alg:bigmerge} or
Algorithm~\ref{alg:assimilate2}.

Let $\hat{w}$ be $w$ with all $s^{\pm 1}$ deleted.   Recall that $D_{r,j} = \Z \cap [j 2^r ,  (j+1) 2^r)$.  For integers $r \geq 0$ and $j \geq 1$, define $v_{r,j}$ to be the subword of $\hat{w}$ whose letters are in positions $D_{r,j}$ in  $\hat{w}$.

By Lemmas~\ref{key cost lemma} and \ref{other key cost lemma}, each call on  Algorithm~\ref{alg:small merge} costs at most
$2176 \, \ell(v_{r,j})^2 = 2176 \, . \, (2^r)^2$ for some $r$ and $j$, where $v_{r,j}$ is the $\theta$  of those lemmas.
 The key point is that this pair  $(r,j)$ never occurs in this way in any other call on Algorithm~\ref{alg:small merge} ---  once that merge has been made it is never repeated.
 So the total cost of all applications of Algorithm~\ref{alg:small merge} is at most $2176 \sum_{\mathcal{D}} 2^{2r}$, where $\mathcal{D}$ is the set of all pairs $(r,j)$ such that $D_{r,j} \subseteq \{ 1, \ldots, \ell(\hat{w}) \}$, and as $\ell(\hat{w}) \leq n$,
\begin{equation*}
 \sum_{\mathcal{D}} 2^{2r} \ \leq \      \sum_{r=0}^{\lceil \log_2 n \rceil } \frac{n}{2^r} 2^{2r}  \ \leq \ n  \sum_{r=0}^{\lceil \log_2 n  \rceil }   2^{r} \ \leq \ n \left( 2^{1+   \lceil \log_2 n  \rceil } -1 \right)  \ \leq \ n(4n-1).
\end{equation*}

Summing our estimates $5n^2 / 2$, $100n^2$,
$15n^2+2n$ and $2176n(4n-1)$, we get the upper bound on cost of $ 17
643n^2 / 2 - 2174n$, which establishes  our theorem.

\bibliography{bibli}
\bibliographystyle{plain}

\end{document}